\documentclass[12pt,twoside]{amsart}
\usepackage{amssymb}

\numberwithin{equation}{section}
\newtheorem{theorem}{Theorem}[section]
\newtheorem{lemma}[theorem]{Lemma}
\newtheorem{proposition}[theorem]{Proposition}
\newtheorem{corollary}[theorem]{Corollary}

\theoremstyle{definition}
\newtheorem{definition}[theorem]{Definition} 
\newtheorem{remark}[theorem]{Remark}
\newtheorem{example}[theorem]{Example}

\begin{document}


\newcommand{\m}[1]{\marginpar{\addtolength{\baselineskip}{-3pt}{\footnotesize \it #1}}}
\newcommand{\A}{\mathcal{A}}
\newcommand{\K}{\mathcal{K}} 
\newcommand{\knd}{\mathcal{K}^{[d]}_n}
\newcommand{\F}{\mathcal{F}}
\newcommand{\N}{\mathbb{N}}
\newcommand{\pr}{\mathbb{P}}
\newcommand{\I}{\mathcal{I}}
\newcommand{\G}{\mathcal{G}}
\newcommand{\lcm}{\operatorname{lcm}}
\newcommand{\ndp}{N_{d,p}}
\newcommand{\tor}{\operatorname{Tor}}
\newcommand{\reg}{\operatorname{reg}} 
\newcommand{\mf}{\mathfrak{m}}
\def\bb{{{\rm \bf b}}}
\def\cc{{{\rm \bf c}}}

 
\title{Sequentially Cohen-Macaulay edge ideals}
\thanks{Version: \today}
\author{Christopher A. Francisco}
\address{Department of Mathematics, University of Missouri, Mathematical Sciences Building, Columbia, MO 65203}
\email{chrisf@math.missouri.edu}
\urladdr{http://www.math.missouri.edu/$\sim$chrisf}
 
\author{Adam Van Tuyl}
\address{Department of Mathematical Sciences \\
Lakehead University \\
Thunder Bay, ON P7B 5E1, Canada}
\email{avantuyl@sleet.lakeheadu.ca}
\urladdr{http://flash.lakeheadu.ca/$\sim$avantuyl/}
 

\begin{abstract}
Let $G$ be a simple undirected graph on $n$ vertices, and let
$\I(G) \subseteq R = k[x_1,\ldots,x_n]$ denote its associated
edge ideal.   We show that all chordal graphs $G$ are sequentially Cohen-Macaulay; 
our proof depends upon showing that the Alexander dual of $\I(G)$ is componentwise linear. 
Our result complements Faridi's theorem that the facet ideal of a simplicial tree is 
sequentially Cohen-Macaulay and implies Herzog, Hibi, and Zheng's theorem that a 
chordal graph is Cohen-Macaulay if and only if its edge ideal is unmixed.  We also characterize 
the sequentially Cohen-Macaulay cycles and produce some examples of 
nonchordal sequentially Cohen-Macaulay graphs.
\end{abstract}
 
\maketitle


\section{Introduction} \label{s.intro}

Let $G$ be a simple graph on $n$ vertices (so $G$ has no loops or multiple edges between two vertices). 
Denote the vertex and edge sets of $G$ by $V_G$ and $E_G$ respectively. We associate to $G$ the 
quadratic squarefree monomial ideal 
$\I(G) \subseteq R = k[x_1,\ldots,x_n]$, with $k$ a field, where 
$\I(G) = (\{x_ix_j ~|~ \{x_i,x_j\} \in E_G \}).$
The ideal $\I(G)$ is called the {\bf edge ideal} of $G$.

The primary focus of this paper is edge ideals of chordal graphs. A graph $G$ is {\bf chordal} if 
every cycle of length $n > 3$ has a chord. Here, if $\{x_1,x_2\},\ldots,\{x_n,x_1\}$ are the $n$ edges 
of a cycle of length $n$, we say the cycle has a chord in $G$ if there exists two vertices 
$x_i,x_j$ in the cycle such that $\{x_i,x_j\}$ is also an edge of $G$, but $\{x_i,x_j\}$ is not an 
edge of the cycle.

We say that a graph $G$ is Cohen-Macaulay if $R/\I(G)$ is Cohen-Macaulay. As 
Herzog, Hibi, and Zheng point out, classifying all the Cohen-Macaulay graphs is probably not tractable 
right now;  this problem is 
as difficult as classifying all Cohen-Macaulay simplicial complexes \cite{HHZ}. However, 
Herzog, Hibi, and Zheng proved in \cite{HHZ} that when $G$ is a chordal graph,
then $G$ is Cohen-Macaulay (over any field) if and only if $\I(G)$ is unmixed.

The property of being sequentially Cohen-Macaulay, a condition
weaker than being Cohen-Macaulay, was introduced by Stanley \cite{Stanley}
 in connection with the theory of nonpure shellability. 

 \begin{definition} \label{d.seqcm}
Let $R=k[x_1,\dots,x_n]$. A graded $R$-module $M$ is called {\bf sequentially Cohen-Macaulay} (over $k$) 
if there exists a finite filtration of graded $R$-modules
\[ 0 = M_0 \subset M_1 \subset \cdots \subset M_r = M \]
such that each $M_i/M_{i-1}$ is Cohen-Macaulay, and the Krull dimensions of the
quotients are increasing:
\[\dim (M_1/M_0) < \dim (M_2/M_1) < \cdots < \dim (M_r/M_{r-1}).\]
\end{definition}

We say that a graph $G$ is sequentially Cohen-Macaulay (over $k$) if $R/\I(G)$ is sequentially 
Cohen-Macaulay. We can expand 
upon Herzog, Hibi, and Zheng's result by using this weakening of the Cohen-Macaulay condition. 
Our main result is the following theorem (which is independent of char$(k)$).

\begin{theorem}[Theorem~\ref{chordalCWL}] All chordal graphs are sequentially Cohen-Macaulay.
\end{theorem}

Thus even chordal graphs whose edge ideals are not unmixed still satisfy a good algebraic property. 
Theorem~\ref{chordalCWL} also generalizes the one-dimensional case of work of Faridi on simplicial 
forests \cite{Faridi}.

Our paper is organized as follows. In the next section, we gather some results from the literature on 
Alexander
 duality and on chordal graphs. In Section~\ref{s.main}, we prove Theorem~\ref{chordalCWL}. We consider 
some nonchordal graphs in Section~\ref{s.cycle}, classifying the sequentially Cohen-Macaulay cycles and 
investigating some properties of graphs containing $n$-cycles for $n > 3$.
We also give a sufficient condition for a graph to fail to be sequentially Cohen-Macaulay.


\section{Required ingredients} \label{s.prelims}

Throughout this paper $G$ will denote a simple graph on $n$ vertices with vertex set $V_G$ and edge 
set $E_G$. Associated to $G$ is the {\bf edge ideal} 
$\I(G) \subseteq R = k[x_1,\ldots,x_n]$, where $\I(G) = (\{x_ix_j ~|~ \{x_i,x_j\} \in E_G \}).$

The {\bf complete graph} on $n$ vertices, denoted $\K_n$, is the graph with edge set
$E_G = \{\{x_i,x_j\} ~|~ 1 \leq i < j \leq n\}$, i.e., the graph with the property that
there is an edge between every pair of vertices.  If $x$ is a vertex of $G$,
we shall write $N(x)$ to denote the neighbors of $x$, that is, those vertices that
share an edge with $x$.  We shall be primarily interested in the case that $G$ is a chordal graph. Chordal graphs have the following property:

\begin{lemma}\cite[Lemma 6.7.12]{V} \label{chordallemma}
Let $G$ be a chordal graph, and let $\K$ be a complete subgraph of $G$.  If $\K \neq G$,
then there is a vertex $x \not \in V_{\K}$ such that the subgraph induced by the 
neighbor set $N(x)$ of $x$ is a complete subgraph.  This also forces the
subgraph induced on $N(x) \cup \{x\}$ to be a complete subgraph.
\end{lemma}

A {\bf vertex cover} of a graph $G$ is a subset $A$ of $V_G$ such that 
every edge of $G$ is incident to at least one vertex of $A$.
Note that we never need to include an isolated vertex in a 
vertex cover. For example, if we have a graph on three vertices 
$x_1$, $x_2$, and $x_3$, and $\{x_1,x_2\}$ is the only edge, then 
$\{x_1\}$ and $\{x_2\}$ are both vertex covers.
The vertex covers of a graph $G$ are related to the Alexander dual
of $\I(G)$.  
\begin{definition}  Let $I$ be a squarefree monomial ideal.
The {\bf squarefree Alexander dual} of 
$I = (x_{1,1}\cdots x_{1,{s_1}},\ldots,x_{t,1}\cdots x_{t,{s_t}})$
is the ideal 
\[I^{\vee} = 
(x_{1,1},\ldots,x_{1,s_1}) \cap \cdots \cap (x_{t,1},\ldots,x_{t,s_t}).\]
\end{definition}

A simple exercise will then verify:

\begin{lemma} 
Let $G$ be a simple graph with edge ideal $\I(G)$.  Then
\[\I(G)^{\vee} = (
\{x_{i_1}\cdots x_{i_k} ~|~ \{x_{i_1},\ldots,x_{i_k}\} ~~\mbox{
is a vertex cover of $G$}\}),\] and the minimal generators of $\I(G)^{\vee}$ correspond to minimal vertex covers.
\end{lemma}

Associated to any homogeneous ideal $I$ of $R$ is a {\bf minimal free graded resolution}
\[ 0 \rightarrow \bigoplus_j R(-j)^{\beta_{h,j}(I)} \rightarrow
\cdots \rightarrow \bigoplus_j R(-j)^{\beta_{1,j}(I)} \rightarrow \bigoplus_j R(-j)^{\beta_{0,j}(I)}
\rightarrow I \rightarrow 0 \]
where $R(-j)$ denotes the $R$-module obtained by shifting the degrees of $R$ by
$j$.
The number $\beta_{i,j}(I)$ is the {\boldmath$ij$}{\bf -th graded Betti number} of $I$ and equals the number of minimal generators of degree $j$ in the $i$-th syzygy module. 
 
\begin{definition} \label{d.linearres}
Suppose $I$ is a homogeneous ideal of $R$ whose generators
all have degree $d$.  Then $I$ has a {\bf linear resolution} if for all $i \ge 0$, $\beta_{i,j}(I) = 0$
for all $j \neq i+d$.
\end{definition}
 
For a homogeneous ideal $I$, we write $(I_d)$ to
denote the ideal generated by all degree $d$ elements of $I$.  Note that $(I_d)$ is different from $I_d$, 
the vector space of all degree $d$ elements of $I$. Herzog and Hibi introduced the
following definition in \cite{HH}.

\begin{definition}
\label{d:CWL}
A homogeneous ideal $I$ is {\bf componentwise linear} if $(I_d)$ has a linear
resolution for all $d$.
\end{definition}

If $I$ is generated by squarefree
monomials, let $I_{[d]}$ denote the ideal generated by the squarefree monomials
of degree $d$ of $I$.  Herzog and Hibi \cite[Proposition 1.5]{HH} showed:

\begin{theorem} \label{Idlinear}
Suppose $I$ is a monomial ideal generated by
squarefree monomials.  Then $I$ is componentwise linear if and only
if $I_{[d]}$ has a linear resolution for all $d$.
\end{theorem}

One can use linear quotients to determine if an ideal has a linear resolution.

\begin{definition}
\label{d:linearquo}
Let $I$ be a monomial ideal of $R$. We say that $I$ has {\bf linear quotients} if for some ordering
$u_1,\dots,u_m$ of the minimal generators of $I$ with
$\deg u_1 \le \deg u_2 \le \cdots \le \deg u_m$ and all $i > 1$,
 $(u_1,\dots,u_{i-1}):(u_i)$ is generated by a subset of $\{x_1,\dots,x_n\}$.
\end{definition}
 
We then require  \cite[Lemma 5.2]{Faridi}:

\begin{lemma}\label{l.linearq}
If $I = (u_1,\ldots,u_m)$ is a monomial ideal of $R = k[x_1,\ldots,x_n]$ that has
linear quotients, and all the $u_i$ have the same degree, then $I$ has a linear resolution.
\end{lemma}

We end this section by applying these ideas to edge ideals.

\begin{lemma}\label{l.vertexcover}
 If $\I(G)$ is the edge ideal of a graph $G$, then
\[\I(G)^{\vee}_{[d]} = ( \{x_{i_1}\cdots x_{i_d} ~|~ \{x_{i_1},\ldots,x_{i_d}\} ~~\mbox{is
a vertex cover of $G$ of size d}\}).\]
\end{lemma}

\begin{proof}
Since $\I(G)^{\vee}$ is generated by the minimal vertex covers, 
any squarefree monomial of degree $d$ in $\I(G)^{\vee}$ corresponds to a set of $d$
vertices which contains a minimal vertex cover, and thus, the $d$ vertices
also form a vertex cover of $G$. 
\end{proof}

\begin{lemma} \label{l.complete}
Let $G=\K_n$ be the complete graph on $n$ vertices. For each $d$, $\I(G)^{\vee}_{[d]}$
has linear quotients.  Thus $\I(G)^{\vee}$ is componentwise linear.
\end{lemma}

\begin{proof}
We show that for each $d$, $\I(G)^{\vee}_{[d]}$ has 
linear quotients and hence a linear resolution, which means that $\I(G)^{\vee}$ is componentwise linear by Theorem \ref{Idlinear}. 

The minimal vertex covers of $\K_n$ are all subsets of $V_{\K_n}$ of size $n-1$.
Hence, by Lemma \ref{l.vertexcover}, $\I(\K_n)^{\vee}_{[d]} = (0)$ if $d < n-1$ or $d>n$ When $d = n$, $\I(\K_n)^{\vee}_{[d]} = (x_1x_2\cdots x_n)$
is a principal ideal. These cases trivially have linear quotients.  It thus suffices to show
that $\I(\K_n)^{\vee}_{[n-1]}$ has linear quotients.

Note that $\I(\K_n)^{\vee}$ is minimally generated by all squarefree monomials of degree $n-1$,
and hence $\I(\K_n)^{\vee}=\I(\K_n)^{\vee}_{[n-1]}$.  Now $\I(\K_n)^{\vee}$
is a squarefree Veronese ideal and thus has a linear resolution \cite{HHcmdiscrete}. Hence $\I(\K_n)^{\vee}_{[n-1]}$ has linear quotients if one orders the monomials in descending lexicographic order.
\end{proof}

\begin{remark} \label{r.generalcomplete}
A statement more general than Lemma~\ref{l.complete} is true. Let $j \le n$, and let $R=k[x_1,\dots,x_n]$.
 We can consider ideals whose components are all possible ideals generated by $j$ of the $n$ variables:
\[ I = (x_1,\dots,x_j) \cap (x_1,\dots,x_{j-1},x_{j+1}) \cap \cdots \cap (x_{n-j+1},\dots,x_n)\]
We can view these ideals as the Alexander duals of either the Stanley-Reisner ideal of a simplicial 
complex with all possible $(j-2)$-faces but no $(j-1)$-faces or as the facet ideal of a simplicial complex 
with all possible $(j-1)$-faces as its facets. $I$ is minimally generated by all squarefree monomials of
 degree $n-j+1$, and hence it is a squarefree Veronese ideal. Thus $I$ has a linear resolution and is 
therefore componentwise linear.
\end{remark}

For our last lemma we show that to determine if $\I(G)^{\vee}$ is
componentwise linear,  we may reduce to the case in which the graph $G$ has no isolated vertices.

\begin{lemma} \label{l.noisolated}
Let $G$ be a simple graph on $n$ vertices with edge ideal $\I(G) \subseteq R=k[x_1,\dots,x_n]$. Let $H$ be the 
graph $G$ with isolated vertices $x_{n+1},\dots,x_m$ added. Assume that $\I(G)^{\vee}$ is componentwise 
linear. 
Then $\I(H)^{\vee} \subseteq k[x_1,\dots,x_m]$ is componentwise linear.
\end{lemma}

\begin{proof}
Note that the edge ideals of $G$ and $H$ have the same minimal generators, though they live in different 
rings. 
Thus $\I(G)^{\vee}$ and $\I(H)^{\vee}$ have the same minimal generators. By \cite[Lemma 2.9]{FVT}, since 
$\I(G)^{\vee}$ is componentwise linear, $\I(H)^{\vee}$ is also.
\end{proof}


\section{Main theorem} \label{s.main}

In this section we prove the main result of this paper. 
Our proof hinges on the following result of 
 Herzog and Hibi \cite{HH} that links the notions of componentwise linearity and sequential Cohen-Macaulayness.
 
\begin{theorem} \label{t.seqcm}
Let $I$ be a squarefree monomial ideal of $R$.
Then $R/I$ is sequentially Cohen-Macaulay if and only if $I^{\vee}$ is componentwise linear.
\end{theorem}

We have arrived at our main result.

\begin{theorem}\label{chordalCWL}
All chordal graphs are sequentially Cohen-Macaulay.
\end{theorem}

\begin{proof}
Let $G$ be a chordal graph. By Theorem \ref{t.seqcm} it suffices to show that $\I(G)^{\vee}$ is componentwise linear. To show
$\I(G)^{\vee}$ is componentwise linear, we have based our proof on 
Faridi's proof of \cite[Theorem 5.4]{Faridi} that the squarefree part of the facet ideal of a 
simplicial forest has linear quotients in each degree. 
By Theorem \ref{Idlinear}, we need to show that $\I(G)^{\vee}_{[d]}$ has a linear
resolution for each $d$. By Lemma \ref{l.linearq}, it suffices to
show that $\I(G)^{\vee}_{[d]}$ has linear quotients for each $d$.

We induct on the number of vertices in the chordal graph. By Lemma~\ref{l.noisolated}, we may assume that 
$G$ has no isolated vertices. Thus the first case to consider is when we have a graph $G$ on two vertices 
connected by an edge. 
In this case $G = \K_2$, so $\I(G)^{\vee}_{[d]}$ has linear quotients for
each $d$ by Lemma \ref{l.complete}.

Suppose now that $G$ is a chordal graph on $n\geq 3$ vertices that has no isolated vertices (so
$G$ has at least two edges). 
If $G = \K_n$, then we are done by Lemma~\ref{l.complete}.  So, we  
may assume that $G$ is not complete. By Lemma \ref{chordallemma} there is a vertex $x \in V_G$ such
that the induced subgraph on $\{x\} \cup N(x)$ is a complete graph. (For example, take $\K$ to be any
edge of $G$, and then $x$ will be some vertex not incident to that edge.)  Write  
$N(x) = \{y_1,\ldots,y_t\}$.  Observe that $G \backslash \{x\}$
and $G \backslash (N(x) \cup \{x\})$ must be chordal. Note that it is possible that 
$G \backslash (N(x) \cup \{x\})$ is an isolated vertex (or vertices); in this case, its edge ideal is the zero ideal.

Now by Lemma \ref{l.vertexcover}, $\I(G)^{\vee}_{[d]}$ is generated by the squarefree monomials that correspond to
the vertex covers of $G$ of size $d$.
Note that any vertex cover $\{x_{i_1},\ldots, x_{i_d}\}$ of $G$ must cover
the complete subgraph $\K_{t+1}$ formed by $\{x,y_1,\ldots,y_t\}$.  So each vertex cover
must contain at least $t$ vertices of $\{x,y_1,\ldots,y_t\}$.

If  $\{x_{i_1},\cdots, x_{i_d}\}$ is a vertex cover of $G$ that contains $x$, then 
$\{x_{i_1},\ldots,x_{i_d}\} \backslash \{x\}$ must be a vertex cover of $G \backslash \{x\}$.
If a vertex cover $\{x_{i_1},\ldots,x_{i_d}\}$ does not contain $\{x\}$, it
must therefore contain $\{y_1,\ldots,y_t\}$.  But then
$\{x_{i_1},\ldots,x_{i_d}\} \backslash \{y_1,\ldots,y_t\}$ must be a vertex cover of 
$G \backslash (N(x) \cup \{x\})$. (In the case when this subgraph is an isolated vertex, since 
there are no edges, the empty set is a vertex cover, as is any subset of vertices.)

Let $H_1 = G\backslash \{x\}$ and $H_2 = G\backslash (N(x) \cup \{x\})$, and let 
$\I(H_1) \subseteq R_1 = k[x_i ~|~ x_i \in V_{G}\backslash \{x\}]$ and 
$\I(H_2) \subseteq R_2 = k[x_i ~|~ x_i \in V_{G} \backslash \{x,y_1,\ldots,y_t\}]$ be their
respective edge ideals.
From the above discussion, it follows that
\[\I(G)^{\vee}_{[d]} = y_1\cdots y_t\I(H_2)^{\vee}_{[d-t]} + x\I(H_1)^{\vee}_{[d-1]}.\]
Here, we are viewing $\I(H_2)^{\vee}_{[d-t]}$ and $\I(H_1)^{\vee}_{[d-1]}$ as ideals of $R$
with the same generators as $\I(H_2)^{\vee}_{[d-t]} \subseteq R_2$ and 
$\I(H_1)^{\vee}_{[d-1]} \subseteq R_1$.

Since $H_1$ and $H_2$ are both chordal with fewer vertices than $G$, by induction,  
\[\I(H_2)^{\vee}_{[d-t]} = (A_1,\ldots,A_a) ~~\text{and}~~ \I(H_1)^{\vee}_{[d-1]}
= (B_1,\ldots,B_b) \]
have linear quotients. We assume that the $A_i$s and $B_i$s have been written
in the correct order for linear quotients. We now show that 
\[\I(G)^{\vee}_{[d]} = (yA_1,\ldots,yA_a, xB_1,\ldots,xB_b), ~~\mbox{with $y =y_1\cdots y_t$,}\]
has linear quotients with respect to this order of the generators.

Since it is clear that $(yA_1,\ldots,yA_{i-1}):(yA_i)$ has linear quotients for $i=2,\ldots,a$,
we need to check that the following ideal has linear quotients:
\[(yA_1,\ldots,yA_a):(xB_1).\]

First note that because $B_1$ corresponds to a vertex cover of $G \backslash \{x\}$, $B_1$
is divisible by at least $t-1$ of $\{y_1,\ldots,y_t\}$.  (To see this, note that $B_1$ covers the 
complete graph $K_t$ formed by the $y_i$s.)  So there exists at most one $y_{\ell}$
such that $y_{\ell} \not| B_1$.

Now suppose there exist monomials $m$ and $p$ and a $j$ such that
\[mxB_1 = pyA_j.\]
We can assume that $mxB_1$ and $pyA_j$ are squarefree.  There are two cases
to consider.
\vspace{.25cm}

\noindent
{\bf Case 1.}  If $y|B_1$, then $B_1 = yB'_1 = y_1\cdots y_tB'_1$.  Since $B_1$
corresponds to a vertex cover of  $G \backslash \{x\}$, $B'_1$ corresponds
to a vertex cover of size $d-t-1$ of $G \backslash (N(x) \cup \{x\})$.
So $B'_1 \in \I(H_2)^{\vee}_{[d-t-1]}$.  Note that if a variable $z|m$, then $z$ must
be a variable of the ring $R_2$; otherwise $mxB_1$ would not be squarefree.
So, for any variable $z$ such that
$z|m$, $zB'_1 \in \I(H_2)^{\vee}_{[d-t]}$, and hence
$zyB'_1 = zB_1 \in (yA_1,\ldots,yA_a)$.  Thus $zxB_1 \in (yA_1,\ldots,yA_a)$
for any $z$ that divides $m$.
\vspace{.25cm}

\noindent
{\bf Case 2.}  Suppose $y \not| B_1$.  By the observation above, there
exists a $y_{\ell} \in \{y_1,\ldots,y_t\}$ such that
$B_1 = y_1\cdots \hat{y}_{\ell} \cdots y_tB'_1$.  Since $mxB_1 = pyA_j$,
and since $y_{\ell}$ divides the right-hand side, we must have
$y_{\ell}|m$.   Note that $y_{\ell}B_1 = yB'_1$ is cover of $G$ of size $d$
with $B'_1$ a cover of $H_2$ of size $d-t$.  Hence $B'_1 \in \I(H_2)^{\vee}_{[d-t]}$.
Thus $y_{\ell}xB_1 \in (yA_1,\ldots,yA_a)$.

The above two cases imply that $(yA_1,\ldots,yA_a):(xB_1)$ has linear quotients.
To finish the proof, we need to check whether
\[(yA_1,\ldots,yA_a,xB_1,\ldots,xB_{i-1}):(xB_i)\]
is generated by a subset of the
variables. If $mxB_i \in (xB_1,\ldots,xB_{i-1})$ for some monomial $m$, then since $(xB_1,\dots,xB_b)$ 
has linear
quotients, there exists a variable $x_i$ that divides $m$ such that $x_ixB_i \in (xB_1,\ldots,xB_{i-1})$.  
If there is a monomial $m$ such that $mxB_i \in (yA_1,\ldots,yA_a)$, the above
argument can be repeated.
\end{proof}

\begin{remark} \label{r.char}
The proof of Theorem~\ref{chordalCWL} shows that chordal graphs are sequentially Cohen-Macaulay 
regardless of the characteristic of $k$ because the linear quotients property is independent of $k$. Faridi \cite{Faridi2} showed that if $I$ is any monomial ideal
that is sequentially Cohen-Macaulay, then the polarization of $I$, a squarefree monomial ideal
associated to $I$, is also sequentially Cohen-Macaulay.  Thus, if $I$ is any monomial
ideal whose polarization is the edge ideal of a chordal graph, $I$ must be sequentially
Cohen-Macaulay.
\end{remark}

Recall that a graph $G$ is a forest if it has no cycles.  A forest, therefore,
is an example of a chordal graph, so we get:

\begin{corollary}
If $G$ is a forest, then $G$ is sequentially Cohen-Macaulay.
\end{corollary}

\begin{remark}
In \cite{Faridi}, Faridi proved that if $\I(\Delta)$ is the facet ideal of simplicial forest
$\Delta$, then $R/\I(\Delta)$ is sequentially Cohen-Macaulay.  When the simplicial
forest has dimension 1, then $\I(\Delta)$ is simply the edge ideal of a forest.  So,
our result can be viewed as a partial generalization of Faridi's result. 
\end{remark}

We close by describing how our Theorem~\ref{chordalCWL} implies Herzog, Hibi, and Zheng's result characterizing Cohen-Macaulay chordal graphs. We begin with a lemma.

\begin{lemma} \label{l.scmcm}
Let $I$ be a squarefree monomial ideal in $R=k[x_1,\dots,x_n]$. Then $R/I$ is Cohen-Macaulay if and only if $R/I$ is sequentially Cohen-Macaulay and $I$ is unmixed.
\end{lemma}

\begin{proof}
When $R/I$ is Cohen-Macaulay, the result is obvious, so assume that $R/I$ is sequentially Cohen-Macaulay and that $I$ is unmixed. Let $I^{\vee}$ be the Alexander dual of $I$. Then by Theorem~\ref{t.seqcm}, $I^{\vee}$ is componentwise linear since $R/I$ is sequentially Cohen-Macaulay. Moreover, since $I$ is unmixed, $I^{\vee}$ is generated in a single degree, meaning that $I^{\vee}$ actually has a linear resolution. By \cite[Theorem 3]{ER}, $R/I^{\vee \vee}=R/I$ is Cohen-Macaulay.
\end{proof}

Herzog, Hibi, and Zheng's result now follows.

\begin{corollary} \label{c.hhz}
A chordal graph is Cohen-Macaulay if and only if its edge ideal is unmixed.
\end{corollary}

\begin{proof}
All chordal graphs are sequentially Cohen-Macaulay, so the corollary is an immediate consequence of Lemma~\ref{l.scmcm}.
\end{proof}


\section{Sequential Cohen-Macaulayness and nonchordal graphs} \label{s.cycle}

In the previous section we showed that if $G$ is a chordal graph, then $R/\I(G)$
is sequentially Cohen-Macaulay.  We now explore the situation in which $G$ is not
chordal.  As we show, $R/\I(G)$ may or may not be sequentially Cohen-Macaulay.

We begin with a classification of the sequentially Cohen-Macaulay $n$-cycles. Villarreal shows in \cite[Corollary 6.3.6]{V} that the only Cohen-Macaulay cycles have three or five vertices. We prove that these are the only sequentially Cohen-Macaulay cycles as well. Note that this does not follow immediately from Villarreal's result because cycles need not be unmixed (in fact, Exercise 6.2.15
of \cite{V} implies an $n$-cycle is unmixed if and only if $n=3,4,5,7$).

\begin{proposition} \label{p.cycles}
Let $G$ be an $n$-cycle for some $n \geq 3$. Then $G$ is sequentially Cohen-Macaulay 
if and only if $n =3$ or $5$.  In fact, when $n=3,5$, the $n$-cycle is Cohen-Macaulay.
\end{proposition}

\begin{proof}
Since a $3$-cycle is chordal, the result for $n=3$ follows from Theorem~\ref{chordalCWL}, and the Cohen-Macaulayness is easy to see. 
When $n=5$, $\I(G)=(x_1x_2,x_2x_3,x_3x_4,x_4x_5,x_1x_5)$, and $k[x_1,\dots,x_5]/\I(G)$ is Gorenstein.

Now suppose $n=2r$ for $r \ge 2$. We have $2r$ edges to cover, and each vertex is incident to exactly two 
edges. Therefore the minimum cardinality of a vertex cover is $r$, and $\{x_1,x_3,\dots,x_{2r-1}\}$ (odd indices)
 and $\{x_2,x_4,\dots,x_{2r}\}$ (even indices) are the two minimal vertex covers. Thus 
$(\I(G)^{\vee}_{r})=(x_1x_3\cdots x_{2r-1},x_2x_4\cdots x_{2r})$, which is a complete intersection of 
monomials of degree $r \ge 2$, and therefore it does not have a linear resolution. Hence $\I(G)^{\vee}$ 
is not componentwise linear, and $G$ is not sequentially Cohen-Macaulay.

Suppose next that $n=2r+1$ for some $r \ge 3$. 
A minimal vertex cover of $G$ consists of alternating vertices plus one additional vertex since alternating
 vertices leaves a single edge uncovered; hence the lowest degree in which $\I(G)^{\vee}$ is generated is 
degree $r+1$. Therefore there are $2r+1$ minimal generators of degree $r+1$, one for each edge that gets 
double-covered when we add an adjacent vertex. Let $J=(\I(G)^{\vee}_{r+1})$. We show that $J$ does not
 have a linear resolution. This implies that $\I(G)^{\vee}$ is not componentwise linear, and hence $G$ 
is not sequentially Cohen-Macaulay. 

To compute the Betti numbers of $J$, we use simplicial homology. Define a squarefree vector to be a 
vector with its entries in $\{0,1\}$. Let $M$ be a monomial ideal, and let 
\[K^\bb(M)=\, \, \mbox{\{squarefree vectors } \cc \in \{0,1\}^{2r+1} \, \mbox{such that } \frac{x^\bb}{x^\cc} \in M\}.\]
This is the upper Koszul simplicial complex of $M$, defined, for example, in \cite{MS}. We can compute the 
$\mathbb N^n$-graded Betti numbers of $M$ with the relation
\[ \beta_{i,\bb}(M) = \dim_k \tilde{H}_{i-1}(K^\bb(M),k) \]
from \cite[Theorem 1.34]{MS}. Summing over all squarefree $\bb$ with degree $j$ gives $\beta_{i,j}(M)$.

We show that $\beta_{2,2r+1}(J) \not = 0$, which proves that $J$ does not have a linear resolution
 when $r \ge 3$. There is a single squarefree vector corresponding to degree $2r+1$, $\bb=(1,\dots,1)$, 
which is associated to the monomial $m=x_1 \cdots x_{2r+1}$. We have a chain complex
\[ \cdots \longrightarrow C_2(K^\bb(J)) \stackrel{\partial_2}{\longrightarrow} C_1(K^\bb(J)) 
\stackrel{\partial_1}{\longrightarrow} C_0(K^\bb(J)) \stackrel{\partial_0}{\longrightarrow} 
C_{-1}(K^\bb(J)) \longrightarrow 0.\]
Below, we shall use the following notation: If $(i_1,\dots,i_n)$ is a vector with entries in $\{0,1\}$ corresponding to a face in our simplicial complex, we shall often write the face as $[x_{j_1},\dots,x_{j_p}]$, where the $j_t$ are exactly the nonzero entries of $(i_1,\dots,i_n)$. For example, the face $(1,0,0,1,0,1)$ is written as $[x_1,x_4,x_6]$. 

Note that the basis of $C_s(K^\bb(J))$ consists of the $s$-dimensional faces $[x_{i_0},\dots,x_{i_s}]$ of $K^\bb(J)$, and
\[\partial_s([x_{i_0},\dots,x_{i_s}]) = 
\sum_{t=0}^s (-1)^t [x_{i_0}, \dots, \hat{x_{i_t}}, \dots, x_{i_s}] .\]
All the faces with which we work have dimension at most two; we orient the faces so that if 
$i_0 < i_1 < i_2$, we traverse $[x_{i_0},x_{i_1}]$ and $[x_{i_1},x_{i_2}]$ in the positive direction and
 $[x_{i_0},x_{i_2}]$ in the negative direction. Similarly, we direct edges so that going from $x_{i_0}$ to
 $x_{i_1}$ is in the positive direction.

To find $\beta_{2,2r+1}(J)$, we need to compute
 $\dim_k \tilde{H}_1(K^\bb(J),k)=\dim_k (\ker \partial_1/ \text{im } \partial_2)$. If we can produce an 
element in $\ker \partial_1$ that is not in $\text{im } \partial_2$, we will have shown that 
$\beta_{2,2r+1}(J)>0$. We shall refer to vertex covers and the corresponding monomials interchangeably below.

Initially, suppose that $2r+1>7$; we handle the case $2r+1=7$ separately. We claim first that 
$m/x_1x_4x_7 \not \in J$. If it were, then there would be a minimal vertex cover $m'$ that divided it. 
But then $x_2x_3x_5x_6x_8x_{2r+1}$ divides $m'$ since $x_1$, $x_4$, and $x_7$ are missing, and $m'$ is a cover. If $2r+1>9$, then to cover the 
remaining $2r-9$ edges not covered, we need at least $r-4$ vertices. This means that 
$\deg m' \ge 6+r-4=r+2$, but all the minimal vertex covers in $J$ have 
degree $r+1$ since $J=(\I(G)^{\vee}_{r+1})$. Also, when $2r+1=9$, the minimal generators of $J$ have degree five, and 
$x_2x_3x_5x_6x_8x_{2r+1}$ is a {\it minimal} vertex cover of degree six and hence is not divisible by an element of $J$. Thus in either case, $m/x_1x_4x_7 \not \in J$.

Next we show that $m/x_1x_4$, $m/x_4x_7$, and $m/x_1x_7$ are in $J$. To prove this, we need to show 
that a minimal vertex cover divides each of these monomials. In the first case, use
 $x_2x_3x_5x_7\cdots x_{2r+1}$; in the second, $x_2x_3x_5x_6x_8x_{10} \cdots x_{2r}$ works; and 
in the last, use $x_2x_4x_6x_8x_{10}\cdots x_{2r}x_{2r+1}$.

Hence $[x_1,x_4]$, $[x_4,x_7]$, and $[x_1,x_7]$ are edges of $K^\bb(J)$, but $[x_1,x_4,x_7]$ is not a 
face of $K^\bb(J)$. Thus $f=[x_1,x_4]+[x_4,x_7]-[x_1,x_7] \in C_1(K^\bb(J))$ is not in the image of 
$\partial_2$. However, 
\[ \partial_1(f)=[x_4]-[x_1]+[x_7]-[x_4]-([x_7]-[x_1])=0.\]
Thus $f$ is in the kernel of $\partial_1$, and $\beta_{2,2r+1}(J) \not = 0$, so $J$ does not have a 
linear resolution. 

When $2r+1=7$, we need a slightly different argument. One can compute that in this case, the Alexander
 dual of $\I(G)$ is
\[ \I(G)^{\vee} = (x_1x_2x_4x_6, x_1x_3x_4x_6, x_1x_3x_5x_6,
 x_1x_3x_5x_7, x_2x_3x_5x_7, x_2x_4x_5x_7, x_2x_4x_6x_7), \]
and it has minimal graded free resolution
\[ 0 \longrightarrow R(-7) \longrightarrow R(-5)^7 \longrightarrow R(-4)^7 
\longrightarrow \I(G)^{\vee} \longrightarrow 0.\]
Because of the second syzygy in degree seven, $\I(G)^{\vee}=(\I(G)^{\vee}_4)$ does not have a 
linear resolution.
Therefore $G$ is not sequentially Cohen-Macaulay.
\end{proof}

\begin{remark} \label{r.cyclechar}
Proposition~\ref{p.cycles} is independent of the characteristic of $k$. Note that if $k$ has 
prime characteristic, the graded Betti numbers of $R/J$ are either the same as in characteristic zero, 
or they go up since the behavior is the same for the dimensions of the homology 
groups we computed. The dimensions of the homology groups in characteristic $p>0$ are either the 
same as in characteristic zero, or they may increase if there is a $p$-torsion part introduced. 
See, for example, the latter part of the discussion of Universal Coefficients in \cite[Chapter 9]{Rotman}. 
Thus we have $\beta_{2,2r+1}(J)>0$ for $r>2$ over all $k$.
\end{remark}

The case of a 5-cycle shows that the converse of Theorem~\ref{chordalCWL} is false. There are many nonchordal sequentially Cohen-Macaulay graphs. We present two simple examples here to demonstrate 
that small changes in a graph that is not sequentially Cohen-Macaulay can give a graph with the property. For further investigation of this idea, see \cite{FH}.

\begin{example} \label{e.4-edge}
Let $G$ be a 4-cycle, and let $H$ be the graph $G$ with a fifth vertex, connected to $G$ by a single edge. 
Thus $\I(G)=(x_1x_2,x_2x_3,x_3x_4,x_1x_4)$, and $\I(H)=(x_1x_2,x_2x_3,x_3x_4,x_1x_4,x_4x_5)$. 
By Proposition~\ref{p.cycles}, $G$ is not sequentially Cohen-Macaulay. The Alexander dual of $\I(H)$ is
\[ \I(H)^{\vee} = (x_1,x_2) \cap (x_2,x_3) \cap (x_3,x_4) \cap (x_1,x_4) \cap (x_4,x_5)
=(x_2x_4,x_1x_3x_5,x_1x_3x_4).\]
It is easy to check that $\I(H)^{\vee}$ is componentwise linear since it has a single generator in 
degree two and regularity three. Hence $H$ is sequentially Cohen-Macaulay.
\end{example}

\begin{example} \label{e.6-triangle}
For a slightly more complicated example, suppose that $G$ is a 6-cycle, and we obtain the graph $H$ by
 adding a seventh vertex and connecting it to two adjacent vertices of $G$. Thus 
\[ \I(H) = (x_1x_2,x_2x_3,x_3x_4,x_4x_5,x_5x_6,x_1x_6,x_1x_7,x_6x_7), \mbox{ and} \]
\[ \I(H)^{\vee} = (x_2x_4x_6x_7,x_1x_3x_5x_7,x_1x_3x_5x_6,x_1x_3x_4x_6,x_1x_2x_4x_6,x_2x_3x_5x_6x_7,x_1x_2x_4x_5x_7).\]

One can check in Macaulay 2 that $\I(H)^{\vee}$ is componentwise linear, so $H$ is sequentially 
Cohen-Macaulay. We remark that tests in Macaulay 2 suggest that adding a triangle in this way to 
a cycle that is not sequentially Cohen-Macaulay may always produce a sequentially Cohen-Macaulay graph. 
\end{example}

We round out this paper with a sufficient condition for a graph to
fail to be sequentially Cohen-Macaulay.  This condition makes use of 
another characterization of sequential Cohen-Macaulayness of quotients by monomial
ideals due to Duval \cite{Duval}. 

Recall that an element $F \in \Delta$, where $\Delta$ is a simplicial complex, is called
a face of $\Delta$.  The dimension of a face $F$ is $\dim F = |F| -1$.  The dimension
of $\Delta$ is then $\dim \Delta = \max_{F \in \Delta} \{\dim F\}$.  We write $\Delta_i$
to denote the subcomplex of $\Delta$ whose maximal faces (the facets) are all
the faces of $\Delta$ of dimension $i$.

\begin{theorem}[{\cite[Theorem 3.3]{Duval}}] \label{Duval}
Let $I$ be a squarefree monomial ideal, and let $\Delta$ be the simplicial complex defined
by $I$ via the Stanley-Reisner correspondence. Let $\Delta_i$ be the pure $i$-dimensional 
subcomplex of $\Delta$. Then $R/I$ is sequentially Cohen-Macaulay if and
only if for every $i$, $-1 \leq i \leq \dim \Delta$, $R/I_{\Delta_i}$ is Cohen-Macaulay.
\end{theorem}

We also need the following definition \cite{Stanley}.

\begin{definition}
If $\Delta$ is a simplicial complex of dim $d-1$, then the {\boldmath$f$}{\bf -vector} 
$f(\Delta)=(f_{-1},f_0,f_1,\ldots,f_{d-1})$, where $f_i$ is the number of faces of dimension $i$
(where $f_{-1} = 1$).  
If 
\[H_{R/(I_{\Delta})}(t) = \frac{h_0 + h_1t + h_2t^2 + \cdots + h_dt^d}{(1-t)^d}\]
is the Hilbert-Poincare series of $R/I_{\Delta}$, then the {\boldmath$h$}{\bf -vector} of $\Delta$
is $h(\Delta) = (h_0,\ldots,h_d)$.
\end{definition}

The {\bf complement} of a simple graph $G$, denoted $G^c$,
is the graph with the same vertex set as $G$, but with edge set $E_{G^c} =\{\{x_i,x_j\}
~|~ \{x_i,x_j\} \not\in E_G\}$, and the {\bf clique-complex} (sometimes called the flag complex) of a simple graph $H$, denoted $\Delta(H)$, is the simplicial complex whose faces are the subsets of vertices on which the induced subgraph of $H$ is a clique.

\begin{theorem}
Let $G$ be a simple graph.  Let $H_2$ be the set of isolated vertices
of $G^c$, and set $H_1 = G^c \backslash H_2$ (so $G^c$ is the disjoint union
of $H_1$ and $H_2$).  If $\#E_{H_1} - \#V_{H_1} + 1 < 0$, then $\I(G)$ is not 
sequentially Cohen-Macaulay. 
\end{theorem}

\begin{proof}
Since $\I(G)$ is a squarefree monomial ideal, $\I(G)$ also corresponds to a simplicial
complex via the Stanley-Reisner correspondence.   In particular,
$\I(G) = I_{\Delta(G^c)}$ where $\Delta(G^c)$ is the clique-complex associated to $G^c$.
Let $\Delta(G^c)_1$ denote
the pure $1$-dimensional subcomplex of $\Delta(G^c)$.  Now $\Delta(G^c)_1$ is simply the $1$-skeleton
of $G^c$, i.e., it is a graph.  Specifically, $\Delta(G^c)_1 = H_1$.  Since $H_1$ is a graph,
the $f$-vector of $H_1$ is 
\[f(H_1) = (1, \#V_{H_1},\#E_{H_1}).\] 
Using the relation between the $f$-vectors and $h$-vectors as given on page 58 of Stanley's book \cite{Stanley}, we have
\[h(H_1) = (1, \#V_{H_1}-2, \#E_{H_1} - \#V_{H_1}+1)\]
If   $\#E_{H_1} - \#V_{H_1} + 1 < 0$, then $h(H_1)$ has negative values.
So $R/I_{\Delta(G^c)_1}$ is not Cohen-Macaulay by \cite[Corollary 3.2]{Stanley} because
the $h$-vector of a Cohen-Macaulay Stanley-Reisner ring must contain only nonnegative values (in fact, must be an $O$-sequence).
Thus, by Theorem \ref{Duval}, $\I(G) = I_{\Delta(G^c)}$ is not sequentially Cohen-Macaulay.
\end{proof}

\begin{example}
The above result gives an alternative justification for why the $4$-cycle is not 
sequentially Cohen-Macaulay.  Since $G^c = \{\{x_1,x_3\},\{x_2,x_4\}\}$, the graph $G^c$ has two edges, but
4 vertices, so $\I(G)$ cannot be sequentially Cohen-Macaulay  since $2 - 4 + 1 < 0$.
\end{example}

\noindent
{\bf Acknowledgments.} We gratefully acknowledge the computer algebra systems CoCoA \cite{Co} and
Macaulay 2 \cite{M2}, which were invaluable in our work on this paper. 
 We would also like to thank Sara Faridi, H. T\`ai H\`a, Jessica Sidman, and the referee for their comments and suggestions. The research of the second author was supported by a grant from NSERC.



\begin{thebibliography}{99}
\bibitem{Co}
CoCoATeam, CoCoA: a system for doing Computations
 in Commutative Algebra, Available at {\tt http://cocoa.dima.unige.it}

\bibitem{Duval} A. M. Duval, Algebraic shifting and sequentially Cohen-Macaulay simplicial complexes. \emph{Electron. J. Combin.} {\bf 3} (1996), no. 1, Research Paper 21, approx. 14 pp. (electronic).

\bibitem{ER} J. Eagon and V. Reiner, Resolutions of Stanley-Reisner rings and Alexander duality. \emph{J. Pure Appl. Algebra} {\bf 130} (1998), no. 3, 265--275.

\bibitem{Faridi} S. Faridi, Simplicial trees are sequentially Cohen-Macaulay. \emph{J. Pure Appl. Algebra} {\bf 190} (2003), 121--136.

\bibitem{Faridi2} S. Faridi,  Monomial ideals via square-free monomial ideals, 
\emph{Lecture Notes in Pure and Applied Mathematics} {\bf 244} (2005) 85--114.

\bibitem{FH} C. A. Francisco and H. T\`ai H\`a, Whiskers and edge ideals, in preparation.

\bibitem{FVT} C. A. Francisco and A. Van Tuyl, Some families of componentwise linear monomial ideals. (2005) Preprint. {\tt math.AC/0508589}


\bibitem{M2}
D.~R. Grayson and M.~E. Stillman, \emph{Macaulay 2, a software system for
  research in algebraic geometry}.
\newblock \verb|http://www.math.uiuc.edu/Macaulay2/|.

\bibitem{HH} J. Herzog and T. Hibi, Componentwise linear ideals.  \emph{Nagoya Math. J.}  {\bf 153} (1999), 141--153.

\bibitem{HHcmdiscrete} J. Herzog and T. Hibi, Cohen-Macaulay polymatroidal ideals. 
(2004) Preprint. {\tt arXiv:math.AC/0409097}

\bibitem{HHZ} J. Herzog, T. Hibi, and X. Zheng, Cohen-Macaulay chordal graphs. (2004) Preprint. {\tt arXiv:math.AC/0407375}

\bibitem{HT} J. Herzog and Y. Takayama, Resolutions by mapping cones. The Roos Festschrift volume, 2.  \emph{Homology Homotopy Appl.}  {\bf 4}  (2002),  no. 2, part 2, 277--294 (electronic). 

\bibitem{MS} E. Miller and B. Sturmfels, {\it Combinatorial commutative algebra}. Springer, 2005.

\bibitem{Rotman} J. Rotman, {\it An introduction to algebraic topology}. Springer-Verlag, 1988.

\bibitem{Stanley} R. P. Stanley, {\it Combinatorics and commutative algebra. Second edition.} 
Progress in Mathematics {\bf 41}. Birkhäuser Boston, Inc., Boston, MA, 1996.

\bibitem{V} R. Villarreal, {\it Monomial Algebras}. Marcel Dekker, 2001.
\end{thebibliography}
\end{document}